\documentclass[a4paper,12pt]{amsart}
\usepackage{amssymb}
\usepackage{amsmath}
\def\i{\lrcorner}

\def\X{X^{\ast}}
\def\Y{Y^{\ast}}

\def\G{\Gamma}

\def\ni{\noindent}

\def\la{\langle}
\def\ra{\rangle}
\def\.{\cdot}

\def\n{\nabla}
\def\l{\lambda}
\def\t{\tilde}
\def\beq{\begin{eqnarray*}}
\def\eeq{\end{eqnarray*}}
\def\te{\theta}
\def\x{\times}
\def\f{\varphi}

\def\o{\omega}

\def\ld{\ldots}

\def\L{\Lambda}

\def\su{\mathfrak{su}}
\def\g{\mathfrak{g}}

\def\r{\end{proof}}

\def\res{\arrowvert}

\def\LL{{\mathcal L}}

\def \RM{\mathbb{R}}
\def \QM{\mathbb{Q}}

\def \ZM{\mathbb{Z}}
\def \CM{\mathbb{C}}

\def \HM{\mathbb{H}}

\newtheorem{ede}{Definition}[section]

\newtheorem{epr}[ede]{Proposition}

\newtheorem{ath}[ede]{Theorem}

\newtheorem{elem}[ede]{Lemma}

\title{Symmetries of contact metric manifolds}
\author{Florin Belgun, Andrei Moroianu and Uwe Semmelmann}
\thanks{The results in the present article were obtained during our
  visit at Erwin Schr{\"o}dinger Institute in Vienna. We express our
  warmest thanks to the Institute for hospitality and support, and to
  Jos{\'e} Figueroa O'Farrill for his invitation.
 The last two authors are members of the {\sl European Differential Geometry 
Endeavour} (EDGE), Research Training Network HPRN-CT-2000-00101, 
supported by The European Human Potential Programme.} 
\address{Florin Belgun\\ Institut f{\"u}r Mathematik \\
 Universit{\"a}t Leipzig\\
Augustusplatz 10-11\\ D-04109 Leipzig, Germany}
\email{belgun@mathematik.uni-leipzig.de}

\address{Andrei Moroianu \\ CMAT\\ {\'E}cole Polytechnique \\ UMR 7640 du CNRS
\\ 91128 Palaiseau \\ France}
\email{am@math.polytechnique.fr}

\address{Uwe Semmelmann\\ Mathematisches Institut \\ Universit{\"a}t
M{\"u}nchen  \\ Theresienstr. 39 \\ D-80333 M{\"u}nchen, Germany}
\email{Uwe.Semmelmann@mathematik.uni-muenchen.de}

\begin{document}

\begin{abstract}
We study the Lie algebra of infinitesimal isometries on compact Sasakian and
K--contact manifolds. On a Sasakian manifold which
is not a space form or 3--Sasakian, every Killing vector field is an
infinitesimal 
automorphism of the Sasakian structure. For a manifold with K--contact
structure, we prove that there exists a Killing vector field of
constant length which is not an infinitesimal 
automorphism of the structure if and only if the manifold is obtained from the
Konishi bundle of a compact pseudo--Riemannian quaternion--K{\"a}hler
manifold after changing the sign of the metric on a maximal negative
distribution. We also prove that non--regular Sasakian manifolds are
not homogeneous and construct examples with cohomogeneity one. Using 
these results we obtain in the last section the classification of all
homogeneous Sasakian manifolds.

\vspace{0.1cm}
\end{abstract}

\maketitle

\section{Introduction}

There are many interesting geometric situations, which are
characterized by the existence of a Killing vector field
with special properties. Examples are Sasakian manifolds, 
K-contact manifolds or total spaces of $S^1$-bundles.
 
In this article we will study the isometry group of a compact
Riemannian manifold admitting a special unit Killing 
vector field.

The starting point is the simple observation that on a
compact Riemannian manifold $(M,\,g)$ a Killing vector field
$\,\xi\,$ induces a decomposition of the Lie algebra
$\,\g:= Lie (Iso(M))$ of the isometry group into 
eigenspaces of the Lie derivative $\,L_\xi$. The zero eigenspace
$\g_0$ can be identified with the space of Killing vector fields
commuting with $\,\xi$.

In the regular case, i.e. where all orbits of $\,\xi\,$ are
closed and have the same length, it is easy to show (under some
technical restrictions) that
$\,\g_0\,$ is spanned by $\,\xi\,$ and 
the Killing vector fields
of the quotient space $\,M/S^1$, cf. Lemma \ref{l2}.

The space $\g_0$ is thus well understood, at least in the regular
case. We then concentrate on the complement of $\g_0$ in $\g$ with
respect to the decomposition mentioned above. Two interesting remarks
can be made (cf. Lemma \ref{per} below): on the one hand, each Killing
vector field from this complement is orthogonal to $\xi$ at every
point of $M$, and on the other hand, the non--zero eigenvalues
for the action of $\LL_\xi$ on $\g$ are determined by those of the
field of endomorphisms $\f$ corresponding to $d\xi$ on $M$. In
particular, if $\f$  has 
no constant non--zero eigenvalues on $M$, then 
$\g=\g_0$. 

In order to study the complement of $\g_0$ in $\g$ it is thus useful
to make further assumptions on $\f$.  
In the remaining part  of the paper we consider the case where $\f$ defines
at every point of $M$ a complex structure on the orthogonal complement
of $\xi$. The manifolds admitting such a structure are called {\it K--contact},
and among them, those satisfying a further integrability condition are
called {\it Sasakian}. 

For these manifolds the Lie algebra $\,\g\,$ splits (as a vector
space) into
only two subspaces: $\,\g = \g_0 \oplus \g_2$. We provide an easy
proof for the fact that a compact Sasakian manifold with
non-trivial $\,\g_2\,$ has to be a space of constant curvature or
a 3-Sasakian manifold. Moreover we describe the decomposition of
$\,\g\,$ on 3-Sasakian manifolds.

In section 4 we characterize K-contact structures carrying
Killing vector fields of constant length in $\g_2$.
It turns out that these manifolds are obtained from $SO(3)$-bundles
over pseudo-Riemannian quaternion K{\"a}hler manifolds by changing
the sign of the metric. In addition we give an example of
K-contact manifold having Killing vector fields of {\it non-constant
length} in $\,\g_2$.

In section 5 we study the isometry group on irregular Sasakian
manifolds and use this to obtain in the last section the
classification of simply connected compact 
homogeneous Sasakian manifolds. Note that this classification was
already obtained in \cite{bg} under the slightly stronger assumption
that the automorphism group of the Sasakian structure (rather than the 
isometry group of the manifold) acts transitively.

\bigskip

\section{Manifolds with unit Killing vector fields}

Throughout this article we identify vectors and 1--forms on Riemannian
manifolds using the metric.

Let $(M^n,g)$ be a compact Riemannian manifold carrying a Killing 
vector field $\xi$. We denote by $G=Iso(M)$ the isometry group of $M$
and by $\g=\mathfrak{iso}(M)$ the Lie algebra of $G$, identified with the Lie
algebra of Killing vector fields on $M$. Since $G$ is compact, the
adjoint representation of $G$ on $\g$ is orthogonal with respect to
some scalar product on $\g$. In particular, the endomorphism
$\LL_{\xi}$ of $\g$ is antisymmetric, so $i\LL_{\xi}$ is a
Hermitian endomorphism of $\g^\CM$. Let $\g^\CM=\g_0^\CM\oplus\g_{\l_1}^\CM
\oplus\ld\oplus\g_{\l_s}^\CM$ be the decomposition of $\g^\CM$ in
eigenspaces for $i\LL_{\xi}$, where $\g_{\l_k}^\CM$ corresponds to the
eigenvalue $\l_k$ and $\g_0^\CM$ to the eigenvalue $0$. This induces a
direct sum decomposition 
\begin{equation}\label{dec}\g=\g_0\oplus\g_{\l_1}
\oplus\ld\oplus\g_{\l_s},
\end{equation}
where $\LL_{\xi}\res\g_0=0$ and
$\LL_{\xi}\circ\LL_{\xi}\res\g_{\l_k}=-\l_k^2$. Note that all $\l_k$'s
can be chosen to be strictly positive real numbers. For such a
choice, (\ref{dec}) will be called the {\it standard decomposition} of $\g$. 

From now on we will always assume that $\xi$ has constant length 1.

\begin{ede} The vector field $\xi$ is called {\em regular} if its flow
  has closed orbits of constant length. It is called {\em
  quasi--regular} if the flow has closed orbits whose lengths have
  jumps. Finally, $\xi$ is called {\em irregular} if its flow has a
  non--closed orbit.
\end{ede}

It is easy to see that $\xi$ is regular or quasi--regular if and only
if there exists 
some $p>0$ such that $\exp_e(p\xi)=e$, where $e$ is the unit of $G$.

We now study in more detail the case where $\xi$ is regular. The flow of $\xi$
defines an isometric circle action on $M$ and the quotient $N~:=M/S^1$
of $M$ by this action carries a Riemannian metric such that the
projection $M\to N$ is a Riemannian submersion with minimal
fibers. Conversely, every connection $\te$ on a principal
$S^1$--bundle $M$ over a Riemannian manifold $N$ induces a
1--parameter family of metrics $g^t$ on
$M$ turning the bundle projection into a Riemannian submersion with
minimal fibers via the following formula:
\begin{equation}\label{gt}
g^t=g_N+t^2\te\otimes\te.
\end{equation}
The metric $g^t$ carries a unit Killing vector field
whose orbits are closed and have constant
length $2\pi t$. 
The following "folkloric" lemma
relates the exterior derivative $d\xi$ and the curvature of the
connection of the $S^1$--bundle $M\to N$.

\begin{elem}\label{l1}
Let $(M^n,g,\xi)$ be as above and denote by $2\pi\ell$ the length of
the orbits of $\xi$. On the $S^1$--principal fibration $M\to N:=M/S^1$
we define a connection whose horizontal distribution is
$\xi^\perp$. Then the curvature form of 
this connection is equal to $id\xi/\ell$. 

Conversely, if $F$ is a 2--form in $H^2(N,\ZM)$, let $M$ be the
$S^1$--bundle over $N$ with Chern class $[F]$ and let $\xi$ denote the
vector field dual to the connection form of a connection on
$M$ whose curvature form is $iF$. Then for each $t>0$, $\xi_t:=\xi/t$
is a unit Killing vector field for the metric $g^t$ defined by
(\ref{gt}), its orbits have length $2\pi t$ and $d\xi_t=d(t\te)=tF$.

\end{elem}

By rescaling the metric on $M$, we may assume that the orbits of the
Killing vector field $\xi$ have constant length $2\pi$ (this amounts
to take $t=1$ in the construction above). The infinitesimal
isometries of $M$ which preserve $\xi$ (i.e. the space $\g_0$) can be
easily described in the following way:

\begin{elem}\label{l2}{\rm (cf. \cite{kak})} With the notations from Lemma
  \ref{l1}, let $A$ be a Killing vector field on
  $M$ commuting with $\xi$. Then there exists a unique Killing vector field
  $X$ on $N$ and a function $f$ on $N$ (unique up to a constant) such that 
\begin{equation}\label{df}
df=-X\i F
\end{equation}
  and
  $A=\X+f\xi$.  Conversely, if $F$ is harmonic and $N$ is simply
  connected, then for every Killing vector field $X$ on $N$ the equation
  (\ref{df}) has a solution $f$ unique up to a constant and $\X+f\xi$
  is a Killing vector field 
  on $M$ commuting with $\xi$.
\end{elem}
{\it Note:} For a vector field $X$ on $N$, $X^*$ denotes the
horizontal lift of $X$ to $M$, i.e. the unique vector field on $M$
orthogonal to $\xi$ which projects onto $X$.
\begin{proof}
Write $A=f\xi+B$ with $\la B,\xi\ra=0$. Since $\LL_\xi$ preserves the
decomposition $TM=\xi\oplus\xi^\perp$, $f$ is actually a function on
$N$ and $B$ is projectable on some $X\in\G(TN)$. The condition $\n
A(\xi,\xi)=0$ is automatically satisfied, $\n A(\Y,\Y)=0$ for every
$Y\in TN$ is equivalent to $X$ being Killing on $N$ and finally $\n
A(\Y,\xi)+\n A(\xi,\Y)=0$ for every
$Y\in TN$ is equivalent to (\ref{df}). 

To prove the converse we only have to check that $d(X\i F)=0$ for
every Killing vector field $X$ on $N$ and harmonic 2--form $F$. This
follows directly from the Hodge decomposition: the harmonic form
$\LL_X F$ ($\LL_X$ commutes with $\Delta$, because $X$ is Killing) is
equal to $d(X\i F)+X\i dF=d(X\i F)$, which is exact, so they both vanish.
\r

The previous Lemma actually shows that there is a canonical exact
sequence

$$0\to \RM\xi\to\g_0\to \mathfrak{iso}(N)\to 0.$$

We now return to the general setting and prove the following simple
but very useful lemma concerning the complement of $\g_0$ in $\g$ with
respect to the standard decomposition described above. 

\begin{elem}\label{per} Let $A$ be a Killing vector field in
  $\g_{\l_1}\oplus\ld\oplus\g_{\l_s}$. Then

1. $A$ is orthogonal to $\xi$ at every point of $M$. 

2. $\LL_\xi A=-A\i d\xi$.
\end{elem}
\begin{proof} 1. It is enough to prove this for $A$ in some fixed
  $\g_{\l_k}$. As before we can write $A=f\xi+B$ with $\la B,\xi\ra=0$, and
  from the invariance w.r.t. $\LL_\xi$ of the decomposition
  $TM=\xi\oplus\xi^\perp$ we deduce that $\xi(\xi(f))=-\l_k^2f$. On the
  other hand $0=\langle\n_\xi A,\xi\rangle =\xi(f)$, so $f=0$ (as $\l_k\ne 0$).

2. Using the first part of the lemma we get 
\begin{equation}
\LL_\xi A=-\LL_A\xi=-d(A\i \xi)-A\i d\xi=-A\i d\xi.
\end{equation}
\r

If we identify $d\xi$ with an endomorphism of $TM$, we see that $A$
has to be an eigenvector of the symmetric endomorphism $d\xi\circ d\xi$ for
the eigenvalue $-\l_k^2$. This shows that the non-zero coefficients of
$X$ in
the above ``Fourier--type'' decomposition are completely determined by
the algebraic behavior of $d\xi$. 

\ni{\bf Example.} If $d\xi\circ d\xi$ has no constant negative eigenvalue on
$M$, then every Killing vector field on $M$ is an automorphism of the
structure $(M,\xi)$ (i.e. $\g=\g_0$).  

\ni{\bf Example.} For the Hopf fibration $(S^{2n+1},\xi)\to\CM P^n$,
$d\xi$ is twice the pull--back of the K{\"a}hler form of the base, thus
showing that $\g$ is actually reduced to $\g_0\oplus\g_2$.

\ni{\bf Example.} Let $N$ be equal to the product of $s$ copies of
$\CM P^1$ with the Fubini--Study metric (i.e. round spheres of radius
1/2). Take $\l_1,\ld,\l_s$ to be distinct positive integers and let $M$ be the
Riemannian manifold induced (via the procedure described above) by the
$S^1$ bundle over $N$ with curvature $i(\l_1\o_1+\ld+\l_s\o_s)$, where
$\o_i$ is the K{\"a}hler form of the i-th factor. Then the Lie algebra
of infinitesimal isometries of $M$ is equal (as vector space) to
$\g=\g_0\oplus\g_{\l_1}\oplus\ld\oplus\g_{\l_n}$, where $\g_0$ is
isomorphic to $\RM \xi$ plus the direct sum of $s$ copies of $\su_2$ and each
$\g_{\l_i}$ is 2--dimensional.

So far we obtained a description of $\,\g_0 \,$ in the regular case
and made some general remarks on the Killing vector fields in the
orthogonal complement of $\,\g_0\,$. In order to get further results
on this complement it is necessary to impose additional assumptions
on $\,\xi$, more precisely on the algebraic behavior of $\,d\xi\,$. 
In the remaining sections we will consider the case where the
skew--symmetric endomorphism corresponding 
to $\,d\xi\,$ defines a complex structure on the orthogonal complement
$\xi^\perp$. This leads to K-contact structures or, with a further
integrability condition to Sasakian structures. 

\bigskip

\section{Sasakian Manifolds}

A Sasakian structure is a special contact structure on a Riemannian 
manifold. These structures were studied in the seventies by the 
Japanese school and, in the last decade, after the work of B{\"a}r
\cite{ba} and Friedrich et al. \cite{bfgk} on manifolds with
Killing spinors and that of Boyer et al. on 3--Sasakian manifolds
\cite{bgm}, they turned out to constitute one of the most important
special geometries, being the odd--dimensional analogues of K{\"a}hler
manifolds. 

\begin{ede}
A vector field $\xi$ on a Riemannian manifold $(M,g)$ is called a {\it
  Sasakian structure} if $\xi$ is a Killing vector field of unit
  length and
\begin{equation}\label{xd}{\nabla_\cdot}\n\xi=\xi\wedge \cdot.\end{equation}
\end{ede}

In particular, if we apply (\ref{xd}) to two arbitrary vectors and
then take the scalar product with $\xi$  
we find that the tensors
$\varphi:=\n{\xi}$ and $\eta:=g(\xi,.)$ are related by
$$\varphi^2=-Id+\eta\otimes\xi.$$ 

It is an easy exercise (see e.g. \cite{ba}) to check that $(M,g)$ is Sasakian
if and only if the metric cone $(\bar M, \bar g)$ defined by $\bar
M=M\x\RM^+$ and $\bar g=dr^2+r^2g$ is K{\"a}hler. Most experts today
actually prefer to take this last statement as the definition of Sasakian
structures, because it is more geometrical and intuitive.
 
\begin{ede}
A triple $\{\xi_1,\xi_2,\xi_3\}$ of Sasakian structures is called a {\it
  3--Sasakian structure} on $M$ if the following conditions are
satisfied: 

1. The frame $(\xi_1,\xi_2,\xi_3)$ is orthonormal;

2. For each permutation $(i,j,k)$ of signature $\delta$, the tensors
$\varphi_i:=\n{\xi_i}$ and $\eta_i:=g(\xi_i,.)$ are related by
$\varphi_j\,\varphi_i=(-1)^{\delta}\,\varphi_k-\eta_j\otimes\xi_i$. 
\end{ede}

Equivalently, $M$ is 3--Sasakian if and only if the cone $\bar M$ is 
hyperk{\"a}hler, and this can be taken to be the definition as before. 

The next lemma gives a sufficient condition for a manifold to be 
3--Sasakian. It was originally proved in \cite{bfgk}.

\begin{elem}\label{fr}If $\xi_1$ and $\xi_2$ are two orthogonal
vector fields defining Sasakian structures,  
then the triple $\{\xi_1,\xi_2,\xi_3:=\n_{\xi_1}\xi_2\}$ is a 3--Sasakian 
structure.
\end{elem}

It turns out that on compact manifolds, the norm of a Killing vector field 
satisfying 
(\ref{xd}) is a ``characteristic function of the sphere''. In other words, 
we have: 

\begin{elem}\label{sph} {\em (see \cite{ta})} Let $\xi$ be a Killing vector field on a
  compact manifold $M$ satisfying (\ref{xd}). If $M$ has non--constant
  sectional curvature, then $\xi$ has constant length (so it defines a
  Sasakian structure on $M$ after a homothetic change of metric). 
\end{elem} 

The following result is a synthesis of several papers studying the
isometries of Sasakian manifolds \cite{ta}, \cite{a}. We will include
here a short proof using the considerations in Section 2.

\begin{ath}\label{sas}Let $M$ be a Sasakian manifold. 
\begin{itemize}
\item If $M$ is neither 3--Sasakian nor a space form, then every
  infinitesimal isometry of $M$ 
  is an infinitesimal automorphism of the Sasakian structure (i.e. 
  $\g=\g_0$).  
\item If $M$ is 3--Sasakian but has non--constant sectional curvature,
  then there is a Lie algebra decomposition $\g=\g'\oplus\su_2$,
  where $\su_2$ is generated by the three Sasakian vector fields and
  $\g'$ consists of the automorphisms of the 3--Sasakian structure
  (i.e. Killing vector fields commuting with the Sasakian vector fields). 
\end{itemize}
\end{ath}

\begin{proof}
1. From Section 2 it
   follows that $\g=\g_0\oplus \g_2$ (since $d\xi\circ d\xi=-4$). Let $X$
   be a non-zero element of $\g_2$ and $Y:=\frac{1}{2}\LL_\xi X$. Then
   $\LL_Y\xi=-\LL_\xi Y=-2\LL^2_\xi X=2X$, so taking the Lie
   derivative with respect to $Y$ in (\ref{xd}) yields that $X$
   satisfies (\ref{xd}), too. From Lemma \ref{sph} we deduce that
   either $M$ is a space form, or $X$ has constant length (and
   consequently defines a Sasakian structure). By Lemma \ref{fr}, this
   implies that $M$ is 3--Sasakian, a contradiction. So $\g_2=\emptyset$.  

2. Let now $M$ be 3--Sasakian. From the first
   part of the above proof, we deduce that for each $i=1,2,3$, $\g$
   can be decomposed as $\g=\g^i_0\oplus\g^i_2$. Now suppose that
   there is some vector field in say $\g^1_2$ which is not a linear
   combination of $\xi_2$ and $\xi_3$. The arguments above show that
   $X$ defines a Sasakian structure, so the cone $\bar M$ admits both a
   hyperk{\"a}hler structure and another K{\"a}hler structure. On the other
   hand, the cone
   is either irreducible or flat (see \cite{ga}). In the last case,
   $M$ has to be a space form. Otherwise, it cannot be locally symmetric 
   (being Ricci--flat), so the existence of another
   K{\"a}hler structure on a hyperk{\"a}hler manifold yields a further
   holonomy reduction, which is impossible from the Berger Holonomy
   Theorem. This proves that for every permutation $\{i,j,k\}$ of
   $\{1,2,3\}$, $\g^i_2=\la\xi_j,\xi_k\ra$. Denote by $\g':=\cap
   \g^i_0$. It is clear that each $\xi_i$ preserves $\g'$ and that
   $\g=\g'\oplus\la\xi_1,\xi_2,\xi_3\ra$. 
\r

\bigskip

\section{K--contact manifolds}

\begin{ede} 1. A {\em contact metric
  structure} on a Riemannian manifold $(M,g)$ is a unit length vector field
  $\xi$ such that the endomorphism $\f$ defined by
  $g(\f\cdot,\cdot):=\frac{1}{2}d\xi(\cdot,\cdot)$ and the 1--form $\eta:=\la\xi,\cdot\ra$
  are related by
\begin{equation}
\label{kc}\f^2=-1+\eta\otimes\xi.
\end{equation}
In other words, $\f$ defines a complex structure on
  the distribution orthogonal to $\xi$. 

2. A contact metric structure $(M,g,\xi,\f)$ is called {\em K--contact}
   if $\xi$ is Killing. 

3. A contact metric 3--structure is an orthonormal frame of contact
metric structures $(\xi_1,\xi_2,\xi_3)$ such that for each permutation
$(i,j,k)$ of signature $\delta$, the tensors  
$\varphi_i:=-\n{\xi_i}$ and $\eta_i:=g(\xi_i,.)$ are related by
$\varphi_i\,\varphi_j=(-1)^{\delta}\,\varphi_k+\eta_j\otimes\xi_i$.
Equivalently, the endomorphisms $\f_i$ satisfy the quaternionic
relations on the distribution orthogonal to $\{\xi_1,\xi_2,\xi_3\}$. 

\end{ede}

It is well--known (and straightforward to prove) that $M$ has a contact metric
structure if and only if the cone $\bar M$ is almost K{\"a}hler. Using
this simple observation we can easily retrieve (with a shorter proof)
the following result of Kashiwada \cite{ka}: a contact
metric 3--structure is necessarily 3--Sasakian! Indeed, the cone $\bar
M$ of a contact metric 3--structure 
$M$ is almost hyperk{\"a}hler and a lemma by Hitchin \cite{hi} (see below)
states that every almost hyperk{\"a}hler manifold is
hyperk{\"a}hler. So $\bar M$ is hyperk{\"a}hler, i.e. $M$ is 3--Sasakian.

We now construct a family of examples of Riemannian manifolds
admitting 3 orthogonal (non--integrable) K--contact metric structures
which will play an important role in this section. This construction
makes essential use of pseudo--Riemannian geometry.  
 
{\bf Fundamental example.} Let $(Q,h)$ be a pseudo--Riemannian
manifold with holonomy $Sp_{p,q}\.Sp_1\subset SO_{4p,4q}$. The 
simplest examples of such manifolds are, as in the Riemannian case,
discrete quotients of the pseudo--Riemannian quaternionic hyperbolic
spaces $\HM 
H^{p,q}:=Sp_{p,q+1}/Sp_{p-1,q+1}\.Sp_1$ (the
case $p=0$, which is excluded in these examples, can be thought of as
being the quaternionic projective space $(\HM P^q,-can)$, with the
negative metric). Other explicit homogeneous examples were constructed
by Alekseevski and 
Cortes (see \cite{ac}). 

The holonomy condition is equivalent to the
existence of a parallel sub--bundle $E$ of $End(Q)$ locally spanned by
three endomorphisms satisfying the quaternionic relations. The Konishi
construction (see \cite{kon}) carries over verbatim to the
pseudo--Riemannian situation. The $SO_3$--principal bundle $S$
associated to $E$ admits a {\it pseudo 3--Sasakian} metric of
signature $(4p,4q+3)$. (The definition of a pseudo 3--Sasakian
structure is the same as that of a 3--Sasakian structure, except that
the metric is pseudo--Riemannian and the 3 Killing vector fields are
time--like). 

Let us now choose an orthogonal decomposition of the
tangent space of $Q$,  $TQ=D_+\oplus D_-$, such that $h$ is positive
on $D_+$ and negative on $D_-$. On $S$ there is an induced orthogonal
decomposition $TS=D_+^*\oplus D_-^*\oplus V$, where $V$ denotes the
vertical distribution. Obviously, this decomposition is invariant
under each of the Sasakian Killing vector fields. We define a
Riemannian metric $g$ by changing the sign of $h$ on $D_-^*\oplus
V$. The Sasakian vector fields, say $\xi_i$, are still Killing vector
fields for $g$, and the endomorphisms associated to
$\frac{1}{2}d\xi_i$ via the new metric are still
complex structures on $\xi_i^\perp$. 

Nevertheless, this change of
metric has definitely altered the integrability of the contact
structures $(\xi_i,\frac{1}{2}d\xi_i)$, and moreover, the three
K--contact structures $\xi_i$ on $(Q,g)$ do not define a contact metric
3--structure, since the corresponding endomorphisms $\f_i$ do not
satisfy the quaternionic relations on the horizontal distribution
(they satisfy the quaternionic relations on $D_-$ and the
anti--quaternionic relations on $D_+$). The Riemannian structures
obtained in this way are called {\it weakly K--contact
  3--structures}.
 
The following result was proved by Hitchin in the Riemannian case
\cite{hi}. The 
proof given in \cite{hi} works without changes in the
pseudo--Riemannian setting.  

\begin{elem}{\em (pseudo--Riemannian Hitchin Lemma)}\label{hitc} A
  pseudo--Riemannian almost hyperk{\"a}hler manifold is hyperk{\"a}hler. 
\end{elem}
(Three almost complex structures $J_i$ on a pseudo--Riemannian
manifold $(M,h)$ define an almost hyperk{\"a}hler structure if the
$J_i$ satisfy the quaternionic relations and the associated K{\"a}hler
forms $h(J_i\cdot,\cdot)$ are closed.)

Consider now a K--contact structure, denoted $\xi_1$ for later
convenience, on a compact manifold 
$M$. The results in Section 2 show that one can
decompose the Lie algebra $\mathfrak{iso}(M)$ as
$\mathfrak{iso}(M)=\g_0\oplus\g_2$, where  
$\LL_{\xi_1}\res_{\g_0}=0$ and $(\LL_{\xi_1}\circ\LL_{\xi_1})\res_{\g_2}=-4$. 

\begin{ath} Let $(M,g,\xi_1,\f_1)$ be a K--contact manifold and let
  $\mathfrak{iso}(M)=\g_0\oplus\g_2$ be the above decomposition of the
  Lie algebra of infinitesimal isometries of $M$. Then $\g_2$ contains
  a Killing 
  vector field of constant length if and only if $M$ admits a weakly K--contact
  3--structure. 
\end{ath}
\begin{proof} The reverse implication is clear from the previous
  example. Conversely, let $\xi_2\in\g_2$ be a Killing vector field of
  unit length, and denote by $\xi_3:=
  \frac{1}{2}\LL_{\xi_1}\xi_2$, which obviously has unit
  length, too. By the definition of $\g_2$ we have 
  $\LL_{\xi_1}\xi_3=-2\xi_2$. The Killing vector field $\zeta:=[\xi_2,\xi_3]$
  can be computed as follows. On one hand, by the Jacobi identity we
  have $\LL_{\xi_1}\zeta=0$, and
  $$-\la\zeta,\xi_1\ra=d\xi_1(\xi_2,\xi_3)=2\la\n_{\xi_2}\xi_1,\xi_3\ra=-2.$$
Finally, for every $Y\perp\xi_1$, we have
\beq\la\zeta,\n_Y\xi_1\ra&=&-\la\n_Y\zeta,\xi_1\ra=\la\n_{\xi_1}\zeta,Y\ra=
\la\n_{\zeta}\xi_1,Y\ra\\ 
&=&\la-\f_1(\zeta),Y\ra=\la\zeta,\f_1(Y)\ra=-\la\zeta,\n_Y\xi_1\ra,
\eeq
thus showing that $\zeta=2\xi_1$. Consequently, for each permutation
$(i,j,k)$ of $(1,2,3)$ of signature $\delta$ we have
$[\xi_i,\xi_j]=2\delta\xi_k$. Now, if we denote by
$\f_i:=-\n\xi_i=\frac{1}{2}d\xi_i$, $\eta_i:=\la\xi_i,\cdot\ra$ and
take the Lie derivative with respect to $\xi_2$ and $\xi_3$ in
(\ref{kc}) we get
\begin{equation}
\label{ac1}\f_1\f_3+\f_3\f_1=\eta_1\otimes\xi_3+\eta_3\otimes\xi_1,
\end{equation}
and
\begin{equation}
\f_1\f_2+\f_2\f_1=\eta_1\otimes\xi_2+\eta_2\otimes\xi_1.
\end{equation}
Taking $\LL_{\xi_2}$ in (\ref{ac1}) and using
(\ref{kc}) yields $\f_2^2=-1+\eta_2\otimes\xi_2$ and similarly
$\f_3^2=-1+\eta_3\otimes\xi_3$. Finally, taking $\LL_{\xi_3}$ in
(\ref{ac1}) gives 
\begin{equation}
\f_2\f_3+\f_3\f_2=\eta_2\otimes\xi_3+\eta_3\otimes\xi_2,
\end{equation}
which shows that the $\f_i$'s are three anti--commuting complex structures
on the distribution $D$ orthogonal to $\xi_1,\xi_2,\xi_3$. From now on
we restrict our attention to the distribution $D$. The endomorphism
$\f_1\f_2\f_3$ has square 1, so we can decompose $D$ as 
$D=D_+\oplus D_-$, where $D_\pm$ is the eigenspace of $\f_1\f_2\f_3$
corresponding to the eigenvalue $\pm1$. That is, the restrictions of
$\f_i$ to $D_+$ (resp. $D_-$) satisfy the anti--quaternionic
(resp. quaternionic) relations. If $D_+$ is empty, we are done, since
in that case $(\xi_1,\xi_2,\xi_3)$ define a contact metric
3--structure, so by the theorem of Kashiwada, $M$ has to be
3--Sasakian. 

The interesting case is when $D_+$ is not empty. We then
define a pseudo--Riemannian metric $h$ on $M$ by changing the sign of $g$
on $D_+$. The three vector fields $\xi_i$ are still Killing vector
fields for $h$ and the endomorphisms associated to $\frac{1}{2}d\xi_i$ via $h$
satisfy the quaternionic relations on $D$. The cone $(\bar M,\bar h)$ is
then a pseudo--Riemannian almost hyperk{\"a}hler manifold. From  Lemma
\ref{hitc} we deduce that  $(\bar M,\bar h)$ is hyperk{\"a}hler, so $(M,h)$
is  pseudo 3--Sasakian. From the construction of $h$ it follows
immediately that $(M,g)$ has a weakly K--contact 3--structure.
\r

In the Sasakian case we have seen that the component $\g_2$ of the Lie
algebra of Killing fields is either $0$-- or $2$--dimensional, or the
manifold has constant sectional curvature. Moreover, in the first
case, any Killing vector field in $\g_2$ has constant length. The
K-contact structures do not have, in general, such a rigidity
property: 

\begin{epr}\label{p0} There exist deformations of the round metric on
  $S^{2n+1}$, $n\ge 3$, which are non Sasakian K-contact structures
  for which the component 
  $\g_2$ of the Lie algebra of the isometry group contains Killing
  fields of non-constant length.
\end{epr}
\begin{proof}
Let $V=V_1\oplus V_2\simeq \RM^{2n+2}$, be decomposed in the
$(2n-2)$-dimensional subspace generated by the first coordinates, and the
4-dimensional subspace generated by the last ones, and let $S^{2n+1}$
be the standard sphere in $V$.

Let $\g\cong\mathfrak{so}(V_1)\oplus\RM\xi$ be the Lie subalgebra of
$\mathfrak{so}(V)$ consisting of the Killing vector fields of
$S^{2n+1}$ depending only on the first $2n-2$ coordinates, plus the
Killing vector field $\xi$, defined by the standard Sasakian structure
of the sphere (and corresponding to the standard complex structure $J$
on $\RM^{2n+2}\cong\CM^{n+1}$). We will deform the metric on $S^{2n+1}$,
keeping it $G$-invariant (where $G$ is the associated Lie subgroup to
$\g$). Because $2n-2\ge 4$, there are elements in $\g$ which do not
commute with $\xi$. On the other hand, these Killing vector fields do
not have constant length, because they all vanish on some totally
geodesic 3-sphere.

We proceed as follows: Let $X_0,J_0 X_0$ be non-trivial vector fields
on $S^3$, seen as the standard sphere in $V_2\cong \RM^4$,
commuting with and orthogonal to $\xi_0$, the Killing field defined by
the standard Sasakian structure on $S^3$ (and corresponding to the
matrix $J_0$, the standard complex structure on $\RM^4$). Such vectors
exist, and are the lifts of some vector fields on $S^2$, the orbit
space of $\xi_0$ on $S^3$. We consider now the vector fields $X,JX$ on
$S^{2n+1}$, depending only on the last 4 coordinates in $\RM^{2n+2}$,
and which project onto $X_0,J_0X_0$ via the projection $V\rightarrow
V_2$. The vectors $X,JX$ are then $G$-invariant. Let $U$ be the open
set of $S^{2n+1}$ where $X$ and $JX$ do not vanish, and let
$F:S^{2n+1}\rightarrow\RM$ be a non-trivial smooth $G$-invariant function
whose support lies in $U$. 

We construct a new metric, $g_F$, on $S^{2n+1}$ defined by 
the field $A$ of symmetric endomorphisms of $TS^{2n+1}$ relating it to
the standard metric. On $U$, $A$ has eigenvalues $e^F,e^{-F}$ and $1$,
corresponding respectively to the eigenvectors $X$, $JX$ and the vectors
orthogonal to them, and on $S^{2m+1}\smallsetminus U$, $A$ is the
identity. Since $A$ is clearly $G$-invariant, $g_F$ is a $G$-invariant
metric. The pair $(g_F,\xi)$ is a K-contact structure because $\xi$ (viewed
both as a vector and as a 1-form) and $d\xi$ (viewed as a 2-form)
remain the same after the above deformation, and the endomorphism
$\phi_F$ of $\xi^\perp$, identified with the form $d\xi$, changes on $X$
and on $JX$, but remains a complex structure on $\xi^\perp$ (namely,
$\phi_F(X)=e^{-2F}JX$ and $\phi_F(JX)=-e^{2F}X$).

Of course, this induces an almost complex structure on $\mathbb{CP}^n$,
the space of orbits of $\xi$, which is never integrable, because it is
not analytic (it coincides with the standard complex structure on the
non-empty interior of the set where $F$ vanishes, and is different on
the open set where $F\ne 0$). So the constructed K-structure is not
Sasakian (this is also a consequence of Theorem \ref{sas}).
\end{proof}

%{\bf Remark. } In dimension 5, this construction is not
%possible:
%\begin{epr}Let $M^5$ be a K-contact manifold for which $\g_2$
%  contains Killing vector fields of non-constant length. Then $M$ is
%  locally isometric to $S^5$, with the standard metric.
%\end{epr}
\bigskip

\section{Irregular Sasakian manifolds}

In this section $(M,g,\xi)$ denotes an irregular Sasakian manifold.
Consider the exponential orbit $O$ of $\xi$ in $G:=Iso(M)$ and denote
by $T$ the closure of $O$ in $G$. Since $G$ is compact, $T$ is a
compact subgroup of $G$. As the closure of an Abelian subgroup, $T$
is itself Abelian, thus isomorphic to a quotient of a torus
$\tilde T$ by a finite subgroup $\G\subset \tilde T$. Let $\L\subset
\g$ be the Lie algebra of $T$ and $\t T$. 

\begin{elem}\label{cen}The subalgebra $\L$ is central in $\g$.
\end{elem}
\begin{proof} As every 3--Sasakian manifold (or Sasakian space form) is at
  least quasi--regular, Theorem \ref{sas} shows that $\g=\g_0$. Thus, for
  every $X\in\g$ and $t\in\RM$, $Ad_{\exp(t\xi)}X=0$. By continuity we
  obtain $Ad_gX=0$  for every $g\in T$, so $ad_YX=0$ for every $Y\in\L$. 
\r

\begin{elem}The closure of the generic orbit on $M$ of the Sasakian
  vector field has dimension larger than 1.  
\end{elem}
\begin{proof} This amounts to say that the generic orbit of the
  compact Lie group $T$ defined above cannot be a circle. Let
  $\zeta\in\Lambda$ be not collinear to $\xi$. Suppose that the
  $T$--orbit of each point in an open set $U$ of $M$ is a circle. Then
  $\zeta$ is a multiple of $\xi$ on $U$, say $\zeta=f\xi$. As $\zeta$
  is Killing as well, we deduce
  that the symmetric part of $df\otimes\xi$ vanishes on $U$, so $\zeta$
  is a constant multiple of $\xi$ on $U$, and hence on $M$, a
  contradiction. 
\r

A classical and non--trivial result by Boothby and Wang \cite{bw}
states that any contact structure whose automorphism group acts
transitively has to be regular. We give here a simple
proof of a slightly stronger version of this result for the case of
Sasakian structures. 

\begin{epr}\label{nh}
An irregular Sasakian manifold is not homogeneous as Riemannian manifold.
\end{epr}

\begin{proof}Let $\zeta\in\L$ be as before and consider a point $x\in
  M$ where $\zeta_x$ is not collinear to $\xi_x$ (the existence of $x$
  follows from the previous lemma). Consider an arbitrary Killing
  vector field $X$. By Lemma \ref{cen} $X$, $\xi$ and $\zeta$ are
  commuting Killing vector fields, hence 
$\la\n_\zeta\xi,X\ra=0$ (trivial application of the Koszul formula). In
  other words, {\it every} Killing vector field is orthogonal at $x$
  to the vector  
$\phi(\zeta)_x=-(\n_\zeta\xi)_x$, which is non-zero by the assumption
  on $\zeta$. Thus 
  $M$ cannot be homogeneous.  
\r

In the remaining part of this section we construct an example of
irregular Sasakian manifold of cohomogeneity one, showing that the
theorem above is optimal. 

{\bf Examples of cohomogeneity one irregular Sasakian manifolds. } 
We use the following equivalent definition for a Sasakian structure:
{\it it is a $K$-contact structure for which the underlying $CR$ structure
is integrable}. Recall that a $CR$ structure is --- in the most
general setting --- a field of complex structures on a field of
hyperplanes, and $J:=\nabla\xi$ on $Q:=\xi^\perp$ gives us such a structure
on a $K$-contact manifold $(M,g,\xi)$. Moreover, the integrability
condition for the $CR$ structure $(M,Q,J)$ is given by
$$N(X,Y)=0,\ \forall X,Y\in Q,\ \mbox{where}\ N:\Lambda^2Q\rightarrow Q,$$
$$4N(X,Y):=[JX,JY]-J[JX,Y]^Q-J[X,JY]^Q-[X,Y],$$
where $X,Y$ are extended to vector fields contained in $Q$, and $A^Q$
denotes the component in $Q$ of the vector field $A$. It is then clear
that $N\equiv 0$ is equivalent to the integrability of the space of
local orbits of $\xi$, which is almost K{\"a}hler for a $K$-contact
structure, and K{\"a}hler for a Sasakian one.

% We are going to change the Sasakian metric on $S^{2n+1}$ as follows:
% we keep the standard $CR$ structure, and change the Reeb vector field
% $\xi_0$, which is regular, with another one, $\xi$, which commutes
% with $\xi_0$, and which is irregular. We adapt the metric such that
% $\xi$ is orthogonal to the hyperplane $Q=Q_0$ and of length $1$, and
% such that it is still $K$-contact. Because of the integrability of the
% standard $CR$ structure on $S^{2n+1}$, the new metric is Sasakian.

Let $\CM^{n+1}\cong V=V_1\oplus V_2$ be a splitting in two complex vector
spaces of dimensions $1$, resp. $n$, and let $J_0$ be the element of
$\mathfrak{u}(n+1)\subset\mathfrak{so}(2n+2)$ corresponding to the
multiplication by $i$ on $V$, and let $J_1$ be the element of
$\mathfrak{u}(n+1)$ corresponding to the multiplication by $i$ on
$V_1$, and acting trivially on $V_2$. They correspond to Killing
vector fields on the round sphere $(S^{2n+1},g_0)\subset \CM^{n+1}$ denoted
by  $T_0$,
resp. $T_1$. 

The first one is the Killing vector field associated to
the standard Sasakian structure on $S^{2n+1}$ and it is easy to see
that $T:=T_0+aT_1$, where $a\in\RM\smallsetminus\QM\cap(0,1)$ is 
a nowhere vanishing Killing vector field, whose generic orbits are
dense in a torus generated by the action of $T_0$ and 
$T_1$. Moreover, $T$ is obviously transverse to $Q:=T_0^\perp$ and
commutes with $T_0$. We define a new metric $g$ on the sphere such
that $T$ has unit length, is orthogonal to $Q$, and $g|_Q:=\alpha
g_0|_Q$, where $\alpha:=g_0(T,T_0)^{-1}$ and $g_0$ is the standard
metric on $S^{2n+1}$. This new metric is $T_0$-
and $T$-invariant, so $\n T$ is a skew-symmetric  endomorphism of $Q$.

If we extend two arbitrary orthogonal vectors $X,Y$
in $Q$ by local vector fields contained in $Q$, orthogonal to each
other  and commuting with
$T_0$, we get by Koszul's formula 
$$2g(\n_X T,Y)=g([X,T],Y)-g([T,Y],X)-g(T,[X,Y])=-g(T,[X,Y]).$$
Here the other terms vanish as $T$ is Killing and $X$ is orthogonal to
$Y$. Similarly we get $2g_0(\n^0_X T_0,Y)=-g_0(T_0,[X,Y])$, and
we can decompose $[X,Y]=g_0(T_0,[X,Y]) T_0+W$, where $W\in Q$. Then
$g(T,[X,Y])=g_0(T_0,[X,Y])g(T,T_0)$. It follows then that $\n_X
T=\n^0_X T_0=J_0(X)$, i.e. it is the standard $CR$ structure on the
sphere, which is integrable. Hence  $g$ is Sasakian.
 
On the other hand, the metric $g$ is invariant under the action of
$U(1)\times U(n)$, where the two factors act on $V_1$, resp. on $V_2$,
and the generic orbits of this group on $S^{2n+1}$ have codimension
1, thus $g$ has cohomogeneity one, as claimed.

% Let $x_1,\dots,x_{2n+2}$ be the coordinates on $\RM^{2n+2}\supset
% S^{2n+1}$, and consider the following vector fields on $\RM^{2n+2}$,
% tangent to the sphere:
% $$\xi_0:=\sum_{k=1}^{n+1}x_{2k-1}\d_{x_{2k}}-x_{2k}\d_{x_{2k-1}},\mbox{
%   and } \xi_1:=x_1\d_{x_2}-x_2\d_{x_1}.$$
% They commute to each other and with the elements of $G\cong SU(n)$
% given by the $J$-invariant orthogonal endomorphisms of the
% $2n$-dimensional space $\{x_1=x_2=0\}$. Thus, the direct product of
% $G$ with $S^1$ (the latter being generated by $\xi$).

% To be continued...

\bigskip

\section{Homogeneous Sasakian Manifolds}

The aim of this section is to prove the following classification result.

\begin{ath}\label{clas} A simply connected compact Sasakian manifold
is homogeneous (as Riemannian manifold) if and only if it is  
the canonical $S^1$--bundle of a homogeneous simply connected compact
Hodge manifold with the metric described in Lemma \ref{l1}. 
\end{ath}
\begin{proof}If $N$ is a homogeneous simply connected compact Hodge manifold, we can assume (after possibly rescaling the metric by an integer constant) that the cohomology class $[F]$ of the K{\"a}hler form is not an integer multiple of any class in $H^2(N,\ZM)$. The Gysin sequence then shows that the canonical $S^1$--bundle $M$ of $N$ is simply connected. It is well--known that $M$ has 
a Sasakian structure, and Lemma \ref{l2} shows that every Killing
vector field $X$ on $N$ induces a Killing vector field $\X$ on $M$,
projectable onto $X$. As the vertical vector field is itself a
(non--vanishing) Killing vector field, it follows that for every point
$x\in M$ and for every vector $A\in T_xM$ there exists a Killing vector field
$\zeta$ on $M$ such that $\zeta_x=A$. Thus the orbit through every point
of the (compact) isometry group of $M$ is open, hence $M$ is
homogeneous. 

Conversely, let $(M,g,\xi)$ be a homogeneous simply connected
compact Sasakian manifold. We can assume that $M$ is not 3--Sasakian,
since in that case $M$ has to be the canonical $S^1$--bundle of a
generalized flag manifold (see \cite{bgm}). From Proposition \ref{nh}, $M$
is either regular or quasi--regular (i.e. the exponential orbit of
$\xi$ in the isometry group is a circle). Denote by $\f_t$ the
isometry induced by $\exp(t\xi)$ on $M$ and let $p$ be the least
positive number such that $\f_p=Id_M$. If $M$ is quasi--regular, there
exists a positive number $q<p$ and a point $x\in M$ such that
$\f_q(x)=x$.  
For every vector $A\in T_xM$ there exists a Killing vector field
$\zeta$ on $M$ such that $\zeta_x=A$. By our assumption, every Killing
vector field on $M$ commutes with $\xi$, so $(\f_q)_*\zeta=\zeta$. This
shows that $\f_q$ has to be the identity on $M$, a contradiction.

We have thus proved that the Sasakian structure is regular. The
quotient $N$ of $M$ by the Sasakian flow is K{\"a}hler (see \cite{bfgk})
and simply connected (by the exact homotopy sequence). From Lemma
\ref{l2} and Theorem \ref{sas} we see that every Killing vector field
on $M$ is of the form $\X+f\xi$, where $X$ is a Killing vector field
on $N$. For every point $x\in N$ and vector $Y\in T_xN$, take some $y$
in the fiber over $x$; as $M$ is homogeneous, there exists a Killing
vector field $\zeta$ on $M$ such that $\zeta_y=\Y$. If we write
$\zeta=\X+f\xi$ with $X$ Killing vector field on $N$, this last
equation shows that $X_x=Y$. So $N$ is homogeneous. 
\r

\ni{\bf Remark.} As Charles Boyer pointed out to us, the above result
was also obtained in \cite{bg} under the slightly stronger assumption
that the automorphism group of the Sasakian structure acts
transitively.

\ni{\bf Remark.} The classification of simply connected
compact homogeneous K{\"a}hler manifolds can be found in \cite{besse}.
Every such manifold has to be an orbit of the
co--adjoint representation of a compact connected Lie group,
endowed with its canonical complex structure. Note
that every such co--adjoint orbit carries several homogeneous K{\"a}hler
metrics, but there is always a canonically defined invariant K{\"a}hler--Einstein
metric, thus at least one Hodge structure. 

\bigskip

 \labelsep .5cm

\end{document}